\magnification=\magstep1
\input amstex
\voffset=-3pc
\hoffset=.25truein

\documentstyle{amsppt}

\vsize=8.5truein
\hsize=6truein

\voffset=-.3pc

\def\bR{\Bbb R}

\def\cA{ A}
\def\cI{ d}

\def\path{\text{path}}

\def\b0{\bold 0}

\def\bR{\Bbb R}

\def\tA{\widetilde A}

\def\b1{\bold 1}

\def\bu{\bold u}
\magnification=\magstep1
\parskip=6pt
\NoBlackBoxes
\topmatter
\title Some Metric and Homotopy Properties of Partial Isometries\endtitle  
\author Lawrence G.~Brown
\endauthor
%\endtopmatter

\abstract{We show that $\|u^*u-v^*v\|\leq \|u-v\|$ for partial isometries $u$ and $v$.
There is a stronger inequality if both $u$ and $v$ are extreme points of the unit ball of a $C^*$--algebra, and both inequalities are sharp.
If $u$ and $v$ are partial isometries in a $C^*$--algebra $A$ such that $\|u-v\|<1$, then $u$ and $v$ are homotopic through partial isometries in $A$.
If both $u$ and $v$ are extremal, then it is sufficient that $\|u-v|| < 2$.
The constants 1 and 2 are both sharp.
We also discuss the continuity points of the map which assigns to each closed range element of $A$ the partial isometry in its canonical polar decomposition.}
\endabstract

\endtopmatter



\bigskip
An element $a$ of a $C^*$--algebra $A$ has closed range if and only if the spectrum of $a^*a$ omits the interval $(0,\epsilon)$ for some $\epsilon>0$.
For any faithful representation $\pi$ of $A$, it is equivalent to say that the operator $\pi(a)$ has closed range.
For $a$ of closed range $\bu (a)$ denotes the partial isometry in its canonical polar decomposition.
It is well known that $\bu(a)$ exists in $A$.

For a unital $C^*$--algebra $A$, $\b1$ denotes its identity element.
And if $A$ is non--unital, $\tA$ denotes the result of adjoining an identity and $\b1$ is the identity of $\tA$.
Also $\b1_H$ denotes the identity operator on the Hilbert space $H$.
If $u$ is a partial isometry in $A$, its left and right defect projections are $\b1-uu^*$ and $\b1-u^*u$, and if $a$ has closed range the defect projections of $a$ are the same as those of $\bu(a)$.
Theorem 1.6.1 of Sakai [S] states that the unit ball of $A$ has extreme points if and only if $A$ is unital, and Theorem 1 of Kadison [K] states that $u$ is an extreme point if and only if $(\b1-uu^*)A (\b1-u^*u)=\{0\}$.
Extreme points are necessarily partial isometries.  The left and right defect ideals of an extreme point are the (closed two-sided) ideals generated by its defect projections.

In [BP1,\S1] Pedersen and the author defined the concept of quasi--invertibility for elements of a unital $C^*$--algebra $A$, and [BP1, Theorem 1.1] gave seven equivalent conditions for an element $a$ to be quasi--invertible.
One of the conditions amounts to the requirement that $a$ have closed range and that $\bu(a)$ be an extreme point of the unit ball of $A$ (in short, $\bu(a)$ is extremal).
Another condition is that there be ideals $I$ and $J$ such that $IJ=\{0\}$, the image of $a$ is $A/I$ is left invertible, and the image of $a$ in $A/J$ is right invertible.

A theme of this paper is to compare results for general partial isometries with stronger results for extremal partial isometries.
This entails also comparing results for general closed range elements with results for quasi--invertible elements.
In particular, by part (ii) of the remark following Theorem 8 below, $a$ is both a continuity point of $\bu$ and an interior point of the domain of $\bu$ if and only if $a$ is quasi--invertible.
This provides additional evidence for the naturality of the quasi--invertibility concept.
(It follows from [BP2, Proposition 7.5] that the domain of $\bu$ has interior points if and only if $A$ is unital.)

\proclaim{Theorem 1}If $u$ and $v$ are partial isometries, then $\|u^*u-v^*v\|\leq \|u-v\|$.
\endproclaim

\demo{Proof}Since every $C^*$--algebra can be embedded in $B(H)$ for some Hilbert space $H$, we may work in $B(H)$.
Let $p=u^*u$ and $q=v^*v$.
Then by Dixmier [Di] and Krein, Krasnosel'ski\u{i}, and Mil'man [KKM] (cf.~also [Da], [H], [P,\S3], and [RS]) we may write $H=H_{00}\oplus H_{01}\oplus H_{10}\oplus H_{11}\oplus H'$ such that all the summands are invariant under $p$ and $q$ and:
\itemitem{(i)}$p_{|H_{00}}=q_{|H_{00}}=0$,
\itemitem{(ii)}$p_{|H_{01}}=0$, $q_{|H_{01}}=\b1_{H_{01}}$,
\itemitem{(iii)}$p_{|H_{10}}=\b1_{H_{10}}$, $q_{|H_{10}}=0$,
\itemitem{(iv)}$p_{|H_{11}}=q_{|H_{11}}=\b1_{H_{11}}$, and 
\itemitem{(v)}$H'$ can be identified with $L^2(S)\oplus L^2(S)$ in such a way that $p_{|H'}$ becomes 
$\pmatrix 1&0\\ 0&0\endpmatrix$ and $q_{|H'}$ becomes $\pmatrix\cos^2\theta&\cos\theta\sin\theta\\ \cos\theta\sin\theta&\sin^2\theta\endpmatrix$.
In (v) $\theta$ is a measurable function from the semifinite measure space $S$ to $(0,\pi/2)$, and the entries of the matrices are multiplication operators on $L^2(S)$.
Some of the five summands of $H$ may be $\{0\}$.

If $H_{01}\neq \{0\}$, then for a unit vector $x$ in $H_{01}$, $\|ux-vx\|=1$, and $\|u-v\|\geq 1\geq \|p-q\|$.
Thus we may assume $H_{01}=H_{10}=\{0\}$, and we also assume $\|u-v\| < \|p-q\|$.
Then $H'\neq \{0\}$, since $\|p-q\| > 0$.
Let $\varphi=\|\theta\|_\infty$ so that $0 < \varphi\leq \pi/2$ and $\|p-q\|=\sin\varphi$.
For $\varphi'<\varphi$, $S_0=\{s\in S\colon\theta(s)>\varphi'\}$ has positive measure.
Let 
$$
\theta'(s)=\cases \theta(s),s\not\in S_0\\ \varphi',s\in S_0\endcases\text{ and let }q'\text{ be the}
$$
approximant to $q$ obtained by replacing $\theta$ with $\theta'$.
There is a partial isometry $w$ such that $w^*w=q'$, $ww^*=q$, and $\|w-q\|\leq 2\sin ((\varphi-\varphi')/2)$.
Then if $v'=vw$ and if $\varphi-\varphi'$ is sufficiently small, we have ${v'}^{*}v'=q'$ and $\|u-v'\| < \|p-q'\|$.

For $f$ and $g$ vectors in Hilbert spaces, we will denote by $f\times g$ the rank one operator $h\mapsto (h,g)f$.
Let $f$ be a unit vector in $L^2(S_0)$, let $e_1=f\oplus 0$ and $e_2=0\oplus f$.
If $r=e_1\times e_1+e_2\times e_2$, then $r$ is a rank two projection which commutes with $p$ and $q'$.
Therefore $ur$ and $v'r$ are both partial isometries, and $\|ur-v'r\|\leq \|u-v'\| < \|p-q'\|=\sin\varphi'=\|pr-q'r\|$.
Note that $ur=g\times e_1$ and $v'r=h\times k$, where $g,h$, and $k$ are unit vectors and $(e_1,k)=\cos\varphi'$.

Thus we are reduced to the following problem:\ \ Let $e_1$ and $k$ be unit vectors such that $(e_1,k)=\cos\varphi'$ for $0<\varphi'< \pi/2$.
Show that the minimum of $\|g\times e_1-h\times k\|$, taken over unit vectors $g$ and $h$, is $\sin\varphi'$.
Note that $k=\cos\varphi' e_1+\sin\varphi' e_2$ where $(e_i,e_j)=\delta_{ij}$, and 
$$
\|g\times e_1-h\times k\|^2=\|(g\times e_1-h\times k)^*(g\times e_1-h\times k)\|=\|t\|.
$$
The operator $t$ is given by the $2\times 2$ matrix 
$$
\pmatrix 1+\cos^2\varphi'-2x\cos\varphi'&\cos\varphi'\sin\varphi'-(x-iy)\sin\varphi'\\
\cos\varphi'\sin\varphi'-(x+iy)\sin\varphi'&\sin^2\varphi'\endpmatrix
$$
where $(g,h)=x+iy$, $x,y\in\bR$.
We are trying to minimize the larger eigenvalue of $t$ subject to the condition $x^2+y^2\leq 1$.
We may assume $y= 0$, since if we change $y$ to 0, we make the determinant larger and keep the trace the same.
This moves the eigenvalues closer to one another.
Now if $y=0$, the larger eigenvalue is $1-x\cos\varphi'+|x-\cos\varphi'|$.
It is easily seen that the minimum occurs for $x=\cos\varphi'$.
\enddemo

\example{Remark}Since every projection is a partial isometry, it is obvious that the inequality in the theorem is sharp.

The following result is probably not new, possibly folklore, but we don't know a reference
\endexample

\proclaim{Proposition 2}Let $f\colon X\to A$ be a continuous function from a topological space into a $C^*$--algebra such that $f(x)$ has closed range for each $x$ in $X$.
Let $u(x)=\bu(f(x))$, and let $x_0$ be in $X$.
Then the following are equivalent:
\itemitem{(i)}The function $u$ is continuous at $x_0$.
\itemitem{(ii)}The function $x\mapsto u(x)^* u(x)$ is continuous at $x_0$.
\itemitem{(iii)}There are $\epsilon > 0$ and a neighborhood $V$ of $x_0$ such that the spectrum of $|f(x)|$ omits $(0,\epsilon)$ for each $x$ in $V$.
\endproclaim

\demo{Proof}(i) $\Rightarrow$ (ii) is obvious.

(ii) $\Rightarrow$ (iii):\ \ There is a neighborhood $W$ of $x_0$ such that $\|u(x)^* u(x)-u(x_0)^* u(x_0)\| < 1$ for $x$ in $W$.
For each such $x$ the partial isometry $r(x)=\bu(u(x)^* u(x) u(x_0)^* u(x_0))$ satisfies $r(x)^* r(x)=u(x_0)^* u(x_0)$ and $r(x) r(x)^*=u(x)^* u(x)$.
Of course $r(x)$ depends only on $u(x)^* u(x)$, and it is well known that $r(x)$ depends continuously on $u(x)^* u(x)$.
Therefore $r$ is continuous at $x_0$.
If $g(x)=f(x) r(x)$, then $g$ is continuous at $x_0$, and $g(x)^* g(x)$ is an invertible element of the unital $C^*$--algebra $u(x_0)^* u(x_0) A u(x_0)^* u(x_0)$.
If $0<\epsilon < \epsilon_1$, and the spectrum (in $A$) of $|f(x_0)|$ omits $(0,\epsilon_1)$, then it follows that the spectrum of $|g(x)|$ omits $(0,\epsilon)$ for $x$ in some neighborhood $V$ of $x_0$, $V\subset W$.
Since $|f(x)|=r(x)|g(x)| r(x)^*$, the result follows.

(iii) $\Rightarrow$ (i):\ \ Let $h\colon [0,\infty)\to [0,\infty)$ be a continuous function such that $h(0)=0$ and $h(t)=t^{-1}$ for $t\geq\epsilon$.
Then $u(x)=f(x) h(|f(x)|)$ for $x$ in $V$.
It follows that $u$ is continuous on $V$.
\enddemo

\proclaim{Corollary 3}With the above hypotheses, the set of continuity points of $u$ is open.
\endproclaim

\proclaim{Theorem 4}If $u$ and $v$ are partial isometries in a $C^*$--algebra $A$ such that $\|u-v\| < 1$, then $u$ and $v$ are homotopic through partial isometries in $A$.
\endproclaim

\demo{Proof}Let $p=u^*u$ and $q=v^*v$.
If $f(t)=(1-t)u+tv$ for $0\leq t\leq 1$, then $\|f(t)-u\| < 1$, $\forall t$.
It follows that $f(t)p$ has closed range and right support projection $p$ for all $t$.
Thus Proposition 2 implies that $u$ is homotopic to $w=\bu (vp)$.
Now since $\|p-q\|<1$ by Theorem 1, there is a continuous function $g$ on $[0,1]$ such that $g(t)$ is a projection, $\|g(t)-q\| < 1$, $\forall t$, $g(0)=p$, and $g(1)=q$.
Thus $qg(t)$ is a closed range element with left and right support projections equal to $q$ and $g(t)$, and hence $vg(t)$ is a closed range element with left and right support projections $vv^*$ and $g(t)$.
So it follows from Proposition 2 that $w$ is homotopic to $v$.
\enddemo

\example{Remark}It is easy to find examples where $\|u-v\|=1$ and $u$ and $v$ are not homotopic.
For a trivial example take $u=0$ and $v\neq 0$.
For less trivial examples note that if two projections are homotopic through partial isometries, then they are also homotopic through projections.
\endexample

\proclaim{Theorem 5}If $u$ and $v$ are extremal partial isometries in a (necessarily unital) $C^*$--algebra $A$ such that $\|u-v\| < 2$, then $u$ and $v$ are homotopic through extremal partial isometries in $A$.
\endproclaim

\demo{Proof}Let $I_1$ and $J_1$ be the left and right defect ideals of $u$, and let $I_2$ and $J_2$ be the left and right defect ideals of $v$.
We first claim that $I_1 J_2=\{0\}=I_2 J_1$.
If, say, $I_1 J_2\neq \{0\}$, let $\pi\colon A\to B(H)$ be an irreducible representation such that $\pi (I_1 J_2)\neq \{0\}$.
Then $\pi(J_1)=\pi(I_2)=\{0\}$, $\pi(u)$ is a proper isometry, and $\pi(v)$ is a proper co--isometry.
Let $f(t)=(1-t) u+tv$ for $0\leq t\leq 1$ and note that $\|f(t)-u\|<1$ for $0\leq t\leq 1/2$ and $\|f(t)-v\|<1$ for $1/2 \leq t\leq 1$.
Thus $\pi(f(t))$ has closed range and right support projection equal to $\b1_H$ for $0\leq t\leq 1/2$.
By applying Proposition 2 to this function, we conclude that the left support projection of $\pi(f(1/2))$ is homotopic to $\pi(uu^*)$.
But since $\|\pi (f(1/2))-\pi(v)\|<1$ and $\pi(v)$ is a co--isometry, the left support projection of $\pi(f(1/2))$ is $\b1_H$, a contradiction.

Now let $I=I_1+I_2$ and $J=J_1+J_2$, and note that $IJ=\{0\}$.
If $\pi_L\colon A\to A/I$ and $\pi_R\colon A\to A/J$ are the quotient maps, then $\pi_L(u)$ and $\pi_L(v)$ are co--isometries, and $\pi_R(u)$ and $\pi_R(v)$ isometries.
By considering the functions $\pi_L\circ f$ and $\pi_R\circ f$ separately on $[0,1/2]$ and $[1/2,1]$, as above, we conclude that $\pi_L (f(t))$ is right invertible and $\pi_R (f(t))$ is left invertible for all $t$.
Also the $\epsilon$ in Proposition 2 (iii) can be taken to be $1-\|u-v\|/2$ for both $\pi_L\circ f$ and $\pi_R\circ f$.
Therefore $f(t)$ is quasi--invertible for all $t$ and Proposition 2, with $\epsilon=1-\|u-v\|/2$, gives the desired homotopy between $u$ and $v$.
\enddemo

\example{Remarks}(i)\ If in Proposition 2 $f(x)$ is quasi--invertible for all $x$, then condition (iii) of the proposition is automatically satisfied.
This follows from [BP1, \S1] even if it isn't explicitly stated there.
In fact the largest $\epsilon$ such that $(0,\epsilon)$ is omitted in the spectrum of $|f(x)|$ was denoted by $m_q (f(x))$ in [BP1], and it was shown that $m_q$ is continuous on the set of quasi--invertibles.
In fact $|m_q(a)-m_q(b)|\leq \|a-b\|$.
A very short proof of the theorem could be given by citing this fact.

(ii)\ It is obvious that $\|u-v\|=2$ does not imply $u$ and $v$ homotopic, since always $\|u-v\|\leq 2$.
The homotopy theory of extremal partial isometries was the starting point of [BP2].
\endexample

\proclaim{Corollary 6}Under the hypotheses of the theorem, $u$ and $v$ have the same defect ideals.
\endproclaim

\proclaim{Theorem 7}If $u$ and $v$ are extremal partial isometries in a $C^*$--algebra such that $\cI=\|u-v\|\leq\sqrt{2}$, then
$\|u^* u-v^*v\|\leq \cI (1-\cI^2/4)^{1/2}$.
\endproclaim

\demo{Proof}By Corollary 6, $u$ and $v$ have the same left and right defect ideals $I$ and $J$.
Since $u^* u-v^* v=\b1-v^* v-(\b1-u^* u)$, and since $\b1-u^* u$, $\b1-v^* v\in J$, we may work in $A/I$.
Thus we may assume $u$ and $v$ are co--isometries, and we may as well work in $B(H)$.
There is a unique number $\theta$ in $[0,\pi/2]$ such that $\cI=2\sin (\theta/2)$, and we are to show that $\|p-q\|\leq\sin\theta$.
It is equivalent to show that $\| (\b1-p) q\|\leq\sin\theta$ and $\|(\b1-q)p\|\leq\sin\theta$.
To show the first inequality, for example, let $g$ be a unit vector in $qH$, and let $f$ in $pH$ be such that $uf=qg=h$.
If $r$ is the rank one projection $h\times h$, then $\|ru-rv\|\leq\cI$, $ru=h\times f$, and $rv=h\times g$.
Since $|(f-g)|^2=2-2\text{ Re }(f,g)\leq 4\sin^2 (\theta/2)$, we conclude that Re $(f,g)\geq\cos\theta$.
Since $\|(\b1-p)g\|\leq \| (\b1-(f\times f))g\|=\|g-(g,f)f\|$, it is easy to calculate that $\|(\b1-p) g\|\leq\sin\theta$.
\enddemo

\example{Remarks}(i) The expression $\cI(1-\cI^2/4)^{1\over 2}$ defines a monotone increasing function on $[0,\sqrt{2}]$ and a decreasing function on $[\sqrt{2},2]$.

(ii)\ To see that the inequality is sharp, start with a one--dimensional Hilbert space $H_1$, a two--dimensional Hilbert space $H_2$, and unit vectors $h\in H_1$, $f,g\in H_2$ such that $(f,g)=\cos\theta$ for $\theta$ in $[0,\pi/2]$.
If $u_0=h\times f$ and $v_0=h\times g$, then $u_0$ and $v_0$ are co--isometries from $H_2$ to $H_1$ such that $\|u_0-v_0\|=2\sin (\theta/2)$ and $\|u_0^* u_0-v_0^* v_0\|=\sin\theta$.
Then for an infinite dimensional space $H'$, let $u'=u_0\oplus\b1_{H'}$ and $v'=v_0\oplus\b1_{H'}$, from $H_2\oplus H'$ to $H_1\oplus H'$.
Finally let $w$ be a unitary from $H_1\oplus H'$ onto $H_2\oplus H'$, and let $u=wu'$ and $v=wv'$.
\endexample

\proclaim{Theorem 8}Let $\cA$ be a $C^*$--algebra, $a$ an element with closed range, and $p_0$ and $q_0$ the left and right defect projections of $a$.
Then $a$ is a continuity point of $\bu$ if and only if 0 is the only closed range element of $p_0\cA q_0$.
\endproclaim

\demo{Proof}If $b$ is a non--zero closed range element of $p_0 \cA q_0$, then $a+tb$ has closed range for each scalar $t$, and it is obvious that $\bu (a+tb)$ does not converge to $\bu(a)$ as $t$ approaches 0.

Conversely assume 0 is the only closed range element of $p_0\cA q_0$, and let $p$ and $q$ be the left and right support projections of $a$.
Let $(a_n)$ be a sequence of closed range elements such that $\|a_n-a\|\to 0$.
For $n$ sufficiently large there is a unique $s_n$ in $\cA$ such that $s_n=q s_n p$, $s_n a_n q=q$, and $p a_n s_n=p$.
Moreover, $(s_n)$ is a convergent sequence.
Let $b_n=(\b1-p_0 a_n s_n) a_n$.
Note that $\b1-p_0 a_n s_n$ is invertible, since $p_0 a_n s_n$ is nilpotent.
Thus $b_n$ has closed range, $b_n$ has the same right support projection as $a_n$, and $\|b_n-a_n\|\to 0$ (since $\|p_0 a_n\|\to 0$).
Thus by Proposition 2, it is sufficient to show $\bu (b_n)\to \bu (a)$.
Since $p_0 b_n q=0$, we see that $p_0 b_n q_0$ has closed range.
Therefore $p_0 b_n q_0=0$.
It follows that the left support projection of $b_n$ is $p$.
Now Proposition 2, applied to $b_n^*$ and $a^*$, implies that  $\bu (b_n)\to \bu(a)$.
\enddemo

\example{Remarks }(i)\ It follows that every quasi--invertible element of $\cA$ is a continuity point of $\bu$.
This was already known from [BP1] even though not explicitly stated.

(ii)\ In [BP2, \S7] an element $a$ of a $C^*$--algebra $\cA$ was said to be of persistently closed range if for some $\epsilon > 0$, $\|b-a\|<\epsilon$ implies $b$ has closed range.
It was shown that persistently closed range elements exist if and only $\cA$ is unital, and that the persistently closed range property is closely related to, but weaker than, quasi--invertibility.
In particular, $a$ has persistently closed range if and only if every element of $p_0\cA q_0$ has closed range.
It follows that $a$ is both a continuity point of $\bu$ and has persistently closed range if and only if $a$ is quasi--invertible.

(iii)\ If $u$ is a partial isometry in $\cA$ such that $a=u|a|$, then $a=u|a|$ can be considered to be a weak sort of polar decomposition of $a$.
The reader may wonder why we restricted the domain of $\bu$ to elements of closed range instead of considering this more general kind of polar decomposition.
One reason is that without introducing special assumptions on $\cA$, it is hard to imagine a natural single--valued and $\cA$--valued extension of $\bu$.
Some related remarks follow without proofs.

(a)\ If $\cA$ is a von Neumann algebra, every element $a$ of $\cA$ has a canonical polar decomposition within $\cA$.
Let $\bu_1 (a)$ be the partial isometry in this decomposition.
Then the only continuity points of $\bu_1$ (for the norm topology) are the quasi--invertibles (and these are also the only continuity points of $\bu$).

(b)\ For general $\cA$ consider the canonical embedding $\iota \colon \cA\hookrightarrow A^{**}$, where $\cA^{**}$ is the enveloping von Neumann algebra of $\cA$.
Let $\bu_2(a)$ be the partial isometry in the canonical polar decomposition of $\iota (a)$ within $\cA^{**}$.
Then $\bu_2(a)\in\cA$ if and only if $a$ has closed range.
If $a$ does not have closed range, it may be that $a=u|a|$ for some partial isometry $u$ in $\cA$ such that $\pi(a)=\pi(u)|\pi(a)|$ is the canonical polar decomposition of $\pi(a)$ for some faithful representation $\pi$; but this will not be the case for all faithful representations $\pi$.

(c)\ The only continuity points of $\bu_2$ are the quasi--invertibles of $\cA$.

\endexample

\bye